\documentclass[11pt]{article}
       \usepackage{amsfonts}
       \usepackage{stmaryrd}
       \usepackage{latexsym,amssymb,mathrsfs,fancyhdr}
       \font\tenmsb=msbm10
       \font\sevenmsb=msbm7
       \font\fivemsb=msbm5
       \catcode`\@=11
       \ifx\amstexloaded@\relax\catcode`\@=\active
       \endinput\else\let\amstexloaded@\relax\fi
       \def\spaces@{\space\space\space\space\space}
       \def\spaces@@{\spaces@\spaces@\spaces@\spaces@\spaces@}
       \def\space@.  {\futurelet\space@\relax}
       \space@.   %
       \def\Err@#1{\errhelp\defaulthelp@\errmessage{AmS-TeX error: #1}}
       \def\relaxnext@{\let\next\relax}
       \def\accentfam@{7}
       \def\noaccents@{\def\accentfam@{0}}
       \def\Cal{\relaxnext@\ifmmode\let\next\Cal@\else
       \def\next{\Err@{Use \string\Cal\space only in math mode}}\fi\next}
       \def\Cal@#1{{\Cal@@{#1}}}
       \def\Cal@@#1{\noaccents@\fam\tw@#1}
       \def\Bbb{\relaxnext@\ifmmode\let\next\Bbb@\else
       \def\next{\Err@{Use \string\Bbb\space only in math mode}}\fi\next}
       \def\Bbb@#1{{\Bbb@@{#1}}}
       \def\Bbb@@#1{\noaccents@\fam\msbfam#1}
       \newfam\msbfam
       \textfont\msbfam=\tenmsb
       \scriptfont\msbfam=\sevenmsb
       \scriptscriptfont\msbfam=\fivemsb
\usepackage{dsfont}
\usepackage[german,english]{babel}
\usepackage{amsmath,amssymb}
\usepackage[square, comma, sort&compress, numbers]{natbib}
\usepackage{cases}
\usepackage{mathrsfs}
\usepackage{amsmath}
\usepackage{amsthm}
\usepackage{amsfonts}
\usepackage{amssymb}
\usepackage{latexsym}
\usepackage{fancyhdr}
\usepackage{extarrows}
\usepackage{geometry}
\newtheorem{thm}{Theorem}[section]
\newtheorem{prop}[thm]{Proposition}
\newtheorem{lem}[thm]{Lemma}
\newtheorem{rem}[thm]{Remark}
\newtheorem{iteration lemma}[thm]{iteration Lemma}
\newtheorem{cor}[thm]{Corollary}
\newtheorem{defn}[thm]{Definition}

\newtheorem{eg}[thm]{Example}

\newtheorem*{acknowledgements*}{ACKNOWLEDGEMENTS}

\begin{document}

\setlength{\columnsep}{5pt}
\title{\bf   One-sided $(b, c)$-inverses in rings }
\author{Yuanyuan Ke$^{a}$\footnote{ E-mail:  keyy086@126.com }, \ Jelena Vi\v{s}nji\'{c}$^{c}$\footnote{ E-mail: jelena.visnjic@medfak.ni.ac.rs},
\ Jianlong Chen$^{a}$\footnote{ Corresponding author. E-mail: jlchen@seu.edu.cn } , \\\\
$^{a}$Department of  Mathematics, Southeast University,  Nanjing 210096,  China \\
$^{c}$Faculty of Medicine, University of Ni\v{s}, Bul. dr Zorana Djindji\'{c}a 81, \\ 18000 Ni\v{s}, Serbia
}
     \date{}

\maketitle
\begin{quote}
{\textbf{Abstract}: \small In this paper we introduce a new generalized inverse in a ring -- one-sided $(b, c)$-inverse, derived as an extension of $(b, c)$-inverse. This inverse also generalizes one-sided inverse along an element, which was recently introduced by H. H. Zhu et al. [H. H. Zhu, J. L. Chen, P. Patr\'{i}cio, Further results on the inverse along an element in semigroups and rings, Linear Multilinear Algebra, 64 (3) (2016) 393-403]. Also, here we present one-sided annihilator $(b, c)$-inverse, which is an extension of the annihilator $(b, c)$-inverse. Necessary and sufficient conditions for the existence of these new generalized inverses are obtained. Furthermore, we investigate conditions for the existence of one-sided $(b, c)$-inverse of a product of three elements and we consider some properties of one-sided $(b, c)$-inverses.



}

\end{quote}

\section{Introduction }
      Throughout the  paper, we assume that $R$ is a ring with identity.  An involution of $R$ is any map $*:R\rightarrow R$ satisfying
    $$
    (a^\ast)^\ast=a,\  (ay)^\ast=y^\ast a^\ast,\  (a+y)^\ast=a^\ast+y^\ast,\  \mbox{for any}~  a, y \in R.
    $$
     A $*$-ring $R$ denotes the ring $R$ with an involution $*$.
     Let $b, c\in R$.
    The concept of the $(b, c)$-inverse as a generalization of the Moore-Penrose inverse, the Drazin inverse, the Chipman's weighted
    inverse and the Bott-Duffin inverse, was for the first time  introduced by M. P. Drazin  in 2012 \cite{D}, in the settings of
    rings.  Recall that an element $a\in R$ is said to be $(b, c)$-invertible if there exists $y\in R$ such that
    $$y\in (bRy)\cap (yRc),\quad\quad yab=b, \quad\quad cay=c.$$
    If such $y$ exists, it is unique and it is called the $(b, c)$-inverse of $a$, denoted by $a^\otimes$.
    For more results on $(b, c)$-inverse we refer the reader to see \cite{Drazin-LMA2014,Drazin-LAA14-Uniqueproof,Ke-Filomat,W-C-C}.

In \cite{D}, M. P. Drazin also introduced the hybrid and annihilator $(b, c)$-inverse of $a$. An element $y\in R$ is called the hybrid $(b, c)$-inverse of $a$ if it satisfies the following equations:
$$
yay=y, \quad yR=bR, \quad y^\circ=c^\circ.
$$
If the condition $yR=bR$ from the above equations is replaced  by $^\circ b=~^\circ y$, then $y$ is the annihilator $(b, c)$-inverse of $a$. Actually, an element $y\in R$ is called the annihilator  $(b, c)$-inverse of $a$ if the following hold:
$$
yay=y, \quad ^\circ b=~^\circ y, \quad y^\circ=c^\circ.
$$
It is shown in \cite{D} that if the hybrid (or annihilator, resp.) $(b, c)$-inverse of $a$ exists, it is unique.

     X. Mary in \cite{Mary2011} introduced a new generalized inverse, called the inverse along element. An element $a\in R$ is said to be invertible along $d\in R$ (or Mary invertible) if there exists $y\in R$ such that
    $$yad=d=day, \quad yR \subseteq dR, \quad Ry\subseteq Rd.$$
    If such $y\in R$ exists, it is unique and it is called the inverse along element $d$ (or Mary inverse). This inverse unify some well-known generalized inverses, such as the group inverse, Drazin inverse and Moore-Penrose inverse.  Also, the inverse along element $d$ is a special case of the $(b, c)$-inverse, for $(b, c)=(d, d)$ \cite[Proposition 6.1]{D}. Several authors also have studied this new outer inverse (see \cite{Benitez-Boasso-1507,Benitez-Boasso-1509,Zhu-Chen-Patricio-LMA15104,
    Zhu-Patricio-Chen-Zhang-LMA15105}).

    As an extension of Mary inverse, H. H. Zhu et al. \cite{Zhu-Chen-Patricio-LMA15104} recently introduced  new inverse in a ring $R$, called  left (right, resp.) inverse along an element.
    Actually, an element $a\in R$ is left (right, resp.) invertible along $d\in R$ (or left (right, resp.) Mary invertible) if there exists $y\in R$ such that
    $$yad=d ~(\mbox{resp.}~ day=d), \quad Ry \subseteq Rd ~(\mbox{resp.}~ yR \subseteq dR).$$


    For the convenience of the reader, some fundamental concepts are given as follows.

      An element $a\in R$ is said to be Moore-Penrose invertible if there exists $y\in R$ which satisfies the following equations:
      $$(1)~aya=a, \quad (2)~yay=y, \quad (3)~(ay)^*=ay, \quad (4)~(ya)^*=ya.$$
      If such $y$ exists, it is unique and is usually denoted by $a^\dag$. The set of all Moore-Penrose invertible elements of $R$ will be denoted by $R^\dag$.

      An element $a\in R$ is (von Neumann) regular if it has an inner inverse $y$, i.e. if there exists $y\in R$ such that $aya=a$. Any inner inverse of $a$ will be denoted by $a^-$. The set of all regular elements of $R$ will be denoted by $R^-$. 
     If $\delta \subseteq \{1, 2, 3, 4\}$ and $y$ satisfies the equations $(i)$ for all $i\in \delta$, then $y$ is an $\delta$-inverse of $a$. The set of all $\delta$-inverse of $a$ is denoted by $a\delta$. Clearly, $a\{1, 2, 3,4\}=\{a^\dag\}$.

Recall that $a\in R$ is left (right, resp.) regular \cite{Azumaya-1954}  if there exists $x\in R$ such that $a=xa^2$ ($a=a^2x$, resp.). If $a$ is both left and right regular, then $a$ is strongly regular. An element $a$ of $R$ is said to be left (right, resp.) $\pi$-regular if there exists $x\in R$ such that $a^n=xa^{n+1}$ ($a^n=a^{n+1}x$, resp.), for some positive integer $n$. An element is strongly $\pi$-regular if it is both left and right $\pi$-regular. An element $a$ is said to be left (right, resp.) $*$-regular if there exists $x\in R$ such that $a=aa^*ax$ ($a=xaa^*a$, resp.). The notions of the Drazin and  group inverse can be referred to the literature \cite{Drazin1958}. It is shown in \cite{Drazin1958} that an element $a\in R$ is Drazin invertible if and only if it is strongly $\pi$-regular. Specially, $a$ is group invertible if and only if it is strongly regular.

     For an element $a \in R$, we define the following image ideals
$$aR=\{ax: x\in R \}, \quad Ra=\{xa: x\in R \},$$
and kernel ideals
$$ a^\circ=\{x\in R: ax=0 \}, \quad ~^\circ a=\{x\in R:xa=0 \}.$$
Let $a, y \in R$. Then $aR = yR$ if and only if there exist $u, v \in R$ such that $a = yu$ and $y = av$. Similarly, $Ra = Ry$ if and only if there exist $s, t \in R$ such that $a = sy$ and $y = ta$.\\
The ring of integers is denoted by $\mathbb{Z}$, and $\mathbb{Z}_n$ stands for the factor ring of $\mathbb{Z}$ modulo $n$, i.e. $\mathbb{Z}_n=\mathbb{Z}/n\mathbb{Z}$, where $n$ is a positive integer.

\medskip

  The results of the paper are organized as follows.
  In Section 2, we define one-sided $(b,c)$-inverses and one-sided annihilator $(b,c)$-inverses. Also, we investigate necessary and sufficient conditions for the existence of these generalized inverses and we state some interesting special cases. In Section 3,  we consider one-sided $(b, c)$-inverses of a product of three elements. Actually, we derive some relations of one-sided $(b,c)$-inverse of $paq$ and one-sided generalized inverse of $pa$ and $aq$, where $a,p,q\in R$. In Section 4, we investigate some properties of one-sided $(b, c)$-inverses, such as the relation between one-sided $(b, c)$-inverses and classical one-sided inverses.


\section{Definition and existence of one-sided $(b, c)$-inverses}
 In this section, we first introduce a class of new generalized inverses in a ring $R$, called a left (right) $(b, c)$-inverse and left (right) annihilator $(b, c)$-inverse. Then we investigate necessary and sufficient conditions for the existence of these generalized inverses.

\begin{defn}\label{ke-onesideinverse-def1}
Let  $b, c\in R$. An element $a\in R$ is said to be left $(b, c)$-invertible if there exists $y\in R$ which satisfies the following equations:
\begin{equation}\label{ke-onesideinverse-eq1}
Ry\subseteq Rc \quad \mbox{and}\quad yab=b.
\end{equation}
In this case $y$ is called a left $(b, c)$-inverse of $a$.

Analogously, an element $a\in R$ is right $(b, c)$-invertible if there exists $y\in R$ such that
\begin{equation}\label{ke-onesideinverse-eq2}
yR\subseteq bR \quad \mbox{and}\quad cay=c.
\end{equation}
In this case $y$ is called a right $(b, c)$-inverse of $a$.
\end{defn}


Let us give some examples of left and right $(b, c)$-inverses.

\begin{eg}\rm{For any ring $R$  (or $*$-ring, resp.) and any  $b, c, y\in R$.
\begin{enumerate}
  \item[(i)] If $a$ is $\{1, 3\}$-invertible, then $a$ is right $(1, a^*)$-invertible. Indeed, 
      by \cite[p. 201]{Hartwig-1976block} we know that $y$ is a $\{1, 3\}$-inverse of $a$ if and only if $a^*ay=a^*$ . Obviously, $Ry\subseteq R$. Hence, $y$ is a right $(1, a^*)$-inverse of $a$.
  \item[(ii)] If $a$ is $\{1, 4\}$-invertible, then $a$ is right $(a^*, 1)$-invertible.
      Analogously as in (i).
  \item[(iii)] If $a$ is invertible along element $d$,
      then $a$ is both left and right $(d, d)$-invertible. If $a$ is $(b, c)$-invertible, then $a$ is both left and right $(b, c)$-invertible.
 \item[(iv)] If $a$ is left (right, resp.) Mary invertible,  then $a$ is left (right, resp.) $(d, d)$-invertible.
\end{enumerate}}
\end{eg}

Notice that the condition $Ry \subseteq Rc$ ($yR \subseteq bR$, resp.) implies $c^\circ \subseteq y^\circ$ ($~^\circ b \subseteq ~^\circ y$, resp.). Thus we give the following definition of a right (left, resp.) annihilator $(b, c)$-inverse of $a$.

\begin{defn}\label{ke-onesideinverse-def3}
 Let $b, c\in R$. An element $a\in R$ is said to be right annihilator $(b, c)$-invertible if there exists $y\in R$ satisfying
\begin{equation}\label{ke-onesideinverse-eq3}
c^\circ \subseteq y^\circ \quad \mbox{and}\quad yab=b.
\end{equation}
In this case $y$ is called a right annihilator $(b, c)$-inverse of $a$.

Similarly, an element $a\in R$ is left annihilator $(b, c)$-invertible if there exists $y\in R$ which satisfies the following equations:
\begin{equation}\label{ke-onesideinverse-eq4}
~^\circ b \subseteq ~^\circ y \quad \mbox{and}\quad cay=c.
\end{equation}
In this case $y$ is called a left annihilator $(b, c)$-inverse of $a$.
\end{defn}


\begin{rem} \rm{Let $a, b, c, y\in R$.
\begin{enumerate}
  \item[(i)] From Definitions \ref{ke-onesideinverse-def1} and \ref{ke-onesideinverse-def3}, obviously we have the following fact: if $a$ is left (right, resp.) $(b, c)$-invertible with a left (right, resp.) $(b, c)$-inverse $y$, then $a$ is right (left, resp.) annihilator $(b, c)$-invertible with a right (left, resp.) annihilator $(b, c)$-inverse $y$.

      However, the converse does not hold in general.
      For example, let $R=\mathbb{Z}$. Then 2 is right annihilator $(0, 2)$-invertible, with right a annihilator $(0, 2)$-inverse 1, but 2 is not left $(0, 2)$-invertible. Indeed, $2^\circ=1^\circ$, but $\mathbb{Z}\nsubseteq 2\mathbb{Z}$.
   \item[(ii)] In general, the condition $c^\circ \subseteq y^\circ$ ($~^\circ b \subseteq ~^\circ y$, resp.) does not imply $Ry\subseteq Rc$ ($yR\subseteq bR$, resp.). For example, let $R=\mathbb{Z}$, $y=1,$ $b=c=2$. Obviously, $c^\circ = y^\circ=b^\circ$, but $\mathbb{Z}\nsubseteq 2\mathbb{Z}$. However, if  $c$ ($b$, resp.) is regular, we have $c^\circ \subseteq y^\circ$ ($~^\circ b \subseteq ~^\circ y$, resp.) if and only if  $Ry\subseteq Rc$ ($yR\subseteq bR$, resp.).
\end{enumerate}}
\end{rem}

From the above remark, we have the following result.

\begin{prop}\label{ke-onesideinverse-prop0.1}
Let  $a, b, c, y\in R$. The following statements hold.
\begin{enumerate}
\item[(i)]If $c\in R^-$, then $a$ is left $(b, c)$-invertible with a left $(b, c)$-inverse $y$ if and only if it is right annihilator $(b, c)$-invertible with a right annihilator $(b, c)$-inverse $y$.
\item[(ii)] If $b\in R^-$, then $a\in R$ is right $(b, c)$-invertible with a right $(b, c)$-inverse $y$ if and only if it is left annihilator $(b, c)$-invertible with a left annihilator $(b, c)$-inverse $y$.
\end{enumerate}
\end{prop}


Next we present necessary and sufficient conditions for existence of a left (right) $(b, c)$-inverse of $a\in R$.

\begin{thm}\label{ke-onesideinverse-thm1}
Let $a, b, c\in R$. Then:
\begin{enumerate}
\item[(i)] $a$ is left $(b, c)$-invertible if and only if $Rb=Rcab$;
\item[(ii)] $a$ is right $(b, c)$-invertible if and only if $cR=cabR$.
\end{enumerate}
\end{thm}
\proof
(i).      Suppose that $a$ is left $(b, c)$-invertible. Then there exists $y\in R$ such that $Ry \subseteq Rc$ and $yab=b$. This implies that $y=rc$ for some $r\in R$. Hence $b=yab=rcab\in Rcab$, which gives $Rb=Rcab.$

      Conversely, if $Rb=Rcab$, there exists $s\in R$ such that
      $b=scab$. Let $y=sc$. Then $b=yab$ and $Ry \subseteq Rc.$ Thus,  $a$ is left $(b, c)$-invertible.

(ii). Similarly as (i).
\qed



\medskip

Note that we can get the existence criterion of a left and right inverse along element \cite[Theorem 2.3 and 2.4]{Zhu-Chen-Patricio-LMA15104} as a direct corollary of Theorem \ref{ke-onesideinverse-thm1}.



\medskip

It is shown in \cite[Theorem 2.2]{D} that $a$ is $(b, c)$-invertible if and only if $Rb=Rcab$ and $cR=cabR$. According to \cite[Proposition 3.3]{W-C-C}, if $a$ is $(b, c)$-invertible, then $b, c\in R^-$. Applying Proposition \ref{ke-onesideinverse-prop0.1} and Theorem \ref{ke-onesideinverse-thm1}, we get the next result.

\begin{cor}\label{ke-onesideinverse-cor3}
Let $a, b, c\in R$. The following statements hold:
\begin{enumerate}
\item[(i)] $a$ is $(b, c)$-invertible if and only if it is both left and right $(b, c)$-invertible;
\item[(ii)] $a$ is $(b, c)$-invertible if and only if $b, c\in R^-$ and $a$ is both left annihilator and right annihilator $(b, c)$-invertible.
\end{enumerate}
\end{cor}




In the next proposition we give the relation between the left $(b, c)$-inverse and the right $(b, c)$-inverse. Since one can easily check  it by using Theorem \ref{ke-onesideinverse-thm1} and Proposition \ref{ke-onesideinverse-prop0.1}, we state it without the proof.

\begin{prop}\label{ke-onesideinverse-prop1.1}
Let $R$ be a $*$-ring  and  $a, b, c\in R$. Then the following statements hold.
\begin{enumerate}
\item[(i)] An element $a$ is left $(b, c)$-invertible if and only if $a^*$ is right $(c^*, b^*)$-invertible;
\item[(ii)] Let $b, c\in R^-$. An element $a$ is left annihilator $(b, c)$-invertible if and only if $a^*$ is right annihilator $(c^*, b^*)$-invertible.
\end{enumerate}
\end{prop}

In the following result we present some consequences of right and left annihilator $(b, c)$-invertibility.

\begin{prop}\label{ke-onesideinverse-prop1}
Let $a, b, c\in R$. Then the following is valid.
\begin{enumerate}
\item[(i)] If $a$ is right annihilator $(b, c)$-invertible, then $b^\circ=(cab)^\circ$;
\item[(ii)] If $a$ is left annihilator $(b, c)$-invertible, then $~^\circ c=~^\circ (cab)$.
\end{enumerate}
\end{prop}

\proof
(i). Suppose that $a$ is right annihilator $(b, c)$-invertible. By Definition \ref{ke-onesideinverse-def3}, there is $y\in R$ such that $c^\circ \subseteq y^\circ$ and $yab=b$. For any $x\in (cab)^\circ$, we have $cabx=0$. Then we obtain $abx\in c^\circ \subseteq y^\circ$, i.e. $yabx=0$, which yields $bx=yabx=0$. This means $(cab)^\circ \subseteq b^\circ$. Thus $b^\circ=(cab)^\circ$.

(ii). Analogously.
\qed

\medskip

In general,  the converse of (i) and (ii) in Proposition \ref{ke-onesideinverse-prop1} does not hold. For example, let $R=\mathbb{Z}$ and let $a=2$, $b=c=1$. Then we have $b^\circ=(cab)^\circ$, but there is no $y\in \mathbb{Z}$ such that $1=b=yab=2y$,  i.e. $2$ is not right annihilator $(1, 1)$-invertible in $\mathbb{Z}$.



Next we give the relation between one-sided $(b, c)$-inverse and the hybrid $(b, c)$-inverse.

\begin{thm}\label{ke-onesideinverse-thm2.1-hybrid}
Let  $a, b, c\in R$.  Then:
\begin{enumerate}
  \item [(i)] $y\in R$ is a right $(b, c)$-inverse and right annihilator $(b, c)$-inverse of $a$ if and only if $y$ is the hybrid $(b, c)$-inverse of $a$;
  \item [(ii)] $y\in R$  is a left $(b, c)$-inverse and left annihilator  $(b, c)$-inverse of $a$ if and only if  $yay=y$,$~^\circ b=~^\circ y,$ $Rc=Ry$.
\end{enumerate}
\end{thm}
\proof
(i). If $y\in R$  is a right and right annihilator $(b, c)$-inverse of $a$, by Definitions \ref{ke-onesideinverse-def1} and \ref{ke-onesideinverse-def3} we have
$$
yR\subseteq bR, \quad cay=c, \quad c^\circ \subseteq y^\circ, \quad yab=b.$$
Thus we have $b=yab\in yR$, and $c^\circ \subseteq y^\circ \subseteq (cay)^\circ=c^\circ$, i.e. $bR=yR$, $c^\circ=y^\circ$. The condition $c=cay$ implies $ay-1\in c^\circ=y^\circ$. Therefore, $y=yay.$

Conversely, if $y$ is the hybrid $(b, c)$-inverse of $a$, we have $yay=y, ~ yR=bR, ~ y^\circ=c^\circ.$ The condition $yR=bR$ gives $b=ys$, for some $s\in R$. Thus, we obtain $b=ys=yays=yab$. Moreover, since $y(ay-1)=0$, we have $ay-1\in y^\circ=c^\circ$. Hence, $c=cay.$

(ii). Similarly as in (i).
\qed




\medskip

Using \cite[Definition 6.2]{D}, Theorem \ref{ke-onesideinverse-thm2.1-hybrid} and  Corollary \ref{ke-onesideinverse-cor3}, we have the following proposition.

\begin{prop}\label{ke-onesideinverse-prop2.2-hybird}
Let $a, b, c\in R$.
\begin{enumerate}
\item[(i)]If $y\in R$ is the annihilator $(b, c)$-inverse of $a$, then $y$ is a right annihilator and left annihilator $(b, c)$-inverse of $a$;
\item[(ii)] If $y\in R$  is a left $(b, c)$-inverse and left annihilator  $(b, c)$-inverse of $a$, then $y$ is the annihilator $(b, c)$-inverse of $a$;
\item[(iii)] If $y\in R$  is a right $(b, c)$-inverse and right annihilator  $(b, c)$-inverse of $a$, then $y$ is the annihilator $(b, c)$-inverse of $a$.
\end{enumerate}
Moreover, if $b, c\in R^-$, then the converse of (i)-(iii) are all valid.
\end{prop}
\proof   (i). Suppose that $y$ is the annihilator $(b, c)$-inverse of $a$. By \cite[Definition 6.2]{D}, we have
     $$
     yay=y, ~^\circ b=~^\circ y,~ y^\circ =c^\circ.
     $$
    Since $yay=y$, we have $(ya-1)y=0$, i.e. $ya-1\in ~^\circ y=~^\circ b$. Therefore, $(ya-1)b=0$ and $yab=b$. Similarly, $y(ay-1)=0$ and $y^\circ =c^\circ$ imply  $cay=c$. Hence, by Definition \ref{ke-onesideinverse-def3}, $a$ is right annihilator and left annihilator $(b, c)$-invertible.

    Moreover, if $b, c\in R^-$, by Corollary \ref{ke-onesideinverse-cor3} (ii), we get that the converse of (i) is valid.

    (ii). Let $y\in R$  be a left $(b, c)$-inverse and left annihilator  $(b, c)$-inverse of $a$. By Theorem \ref{ke-onesideinverse-thm2.1-hybrid} (ii), we have $yay=y$,$~^\circ b=~^\circ y,$ $Rc=Ry$. Since the condition $Rc=Ry$ implies $y^\circ=c^\circ$, we have that $y$ is the annihilator $(b, c)$-inverse of $a$.

    Moreover, if $b, c\in R^-$, then $~^\circ b=~^\circ y$ and $y^\circ=c^\circ$ are equivalent to $~^\circ b=~^\circ y$ and $Rc=Ry$ under the condition $yay=y$. Hence, the converse of (ii) holds.

    (iii). Analogously as (ii).
\qed

\begin{rem}\rm{ Let $a, b, c, y\in R$.
\begin{enumerate}
\item[(i)]
Note that in general the converse of (i), (ii) and (iii) from Proposition \ref{ke-onesideinverse-prop2.2-hybird} doesn't hold. For example, let $R=\mathbb{Z}$. We have that 1 is the annihilator (1, 2)-inverse of 1 and 1 is a left annihilator (1, 2)-inverse of 1, but  1 is not left (1, 2)-inverse of 1. Indeed, $a=b=y=1,$ $c=2$, we have $^\circ b=~^\circ y=y^\circ=c^\circ=0$, but $\mathbb{Z} \nsubseteq 2\mathbb{Z}$.
\item[(ii)] If $a$ is both left $(b, c)$-invertible and right annihilator $(b, c)$-invertible, with  a left $(b, c)$-inverse $y$ and  a right annihilator $(b, c)$-inverse $z$, in general  $y$ doesn't have to be equal to $z$. For example, let $R=\mathbb{Z}_8=\mathbb{Z}/8\mathbb{Z}$,  $a=\overline{5}, b=\overline{0}, c=\overline{2}$. Then $\overline{5}$ is left $(\overline{0}, \overline{2})$-invertible and it is also right annihilator $(\overline{0}, \overline{2})$-invertible. And $y=\overline{4}$ is a left $(\overline{0}, \overline{2})$-inverse of $\overline{5}$, $z=\overline{6}$ is a right annihilator $(\overline{0}, \overline{2})$-inverse of $\overline{5}$. But $\overline{4}\neq\overline{6}$ in $\mathbb{Z}_8$. 
\end{enumerate}
}
\end{rem}

Next we give another existence criterion of a left $(b, c)$-inverse of $a$.

\begin{thm}\label{ke-onesideinverse-thm3}
Let $a, b, c\in R$. Then the following are equivalent:
\begin{enumerate}
\item[(i)] $a$ is left $(b, c)$-invertible;
\item[(ii)] $a^\circ \cap bR=\{0\}$, $abR\cap c^\circ=\{0\}$ and $R=Rca+~^\circ b$;
\item[(iii)] $a^\circ \cap bR=\{0\}$ and $R=Rca+~^\circ b$;
\item[(iv)] $abR\cap c^\circ=\{0\}$ and $R=Rca+~^\circ b$;
\item[(v)] $R=Rca+~^\circ b$.
\end{enumerate}
\end{thm}
\proof (i) $ \Rightarrow$ (ii). If $a$ is left $(b, c)$-invertible, by Definition \ref{ke-onesideinverse-def1}, there is $y\in R$ such that $Ry\subseteq Rc$ and $yab=b$. For any $x\in a^\circ \cap bR$, we have $ax=0$ and $x=bs$, for some $s\in R$. Thus $x=bs=yabs=yax=0$, i.e. $a^\circ \cap bR=\{0\}$.

For any $z\in abR\cap c^\circ$, there is $s'\in R$ such that $z=abs'$ and $cz=0$. Then $z=abs'=a(rcab)s'=arcz=0$, i.e. $abR\cap c^\circ=\{0\}$.

Using Theorem \ref{ke-onesideinverse-thm1} (i), we have $Rb=Rcab$. Thus, there is $r\in R$ such that $b=rcab$. Let  $u=1-rca$.  Then $u\in ~^\circ b$. For any $t\in R$, we have $t=t(rca+u)=trca+tu\in Rca+~^\circ b$, i.e. $R=Rca+~^\circ b$.

(ii) $ \Rightarrow$ (iii)(or (iv), resp.). Obviously.

(iii)(or (iv), resp.)  $ \Rightarrow$ (v). Clearly.

(v) $ \Rightarrow$ (i). The condition $R=Rca+~^\circ b$ gives $b\in Rcab$, so $Rb=Rcab$. According to Theorem \ref{ke-onesideinverse-thm1} (i), $a$ is left $(b, c)$-invertible.
\qed

\medskip

Dually, we have the following result for the existence of a right $(b, c)$-inverse of $a$.

\begin{thm}\label{ke-onesideinverse-thm4}
Let $a, b, c\in R$. Then the following are equivalent:
\begin{enumerate}
\item[(i)] $a$ is right $(b, c)$-invertible;
\item[(ii)] $~^\circ a\cap Rc=\{0\}$, $Rca\cap ~^\circ b=\{0\}$ and $R=abR+ c^\circ$;
\item[(iii)] $~^\circ a\cap Rc=\{0\}$ and $R=abR+ c^\circ$;
\item[(iv)] $Rca\cap ~^\circ b=\{0\}$ and $R=abR+ c^\circ$;
\item[(v)] $R=abR+ c^\circ$.
\end{enumerate}
\end{thm}

Applying Theorem \ref{ke-onesideinverse-thm3}, Theorem \ref{ke-onesideinverse-thm4} and Corollary \ref{ke-onesideinverse-cor3}, we get the following corollary.

\begin{cor}
Let $a, b, c\in R$.  Then the following are equivalent:
\begin{enumerate}
\item[(i)] $a$ is  $(b, c)$-invertible;
\item[(ii)] $a^\circ \cap bR=\{0\}$, $~^\circ a\cap Rc=\{0\}$, $R=abR\oplus c^\circ$ and $R=Rca\oplus ~^\circ b$;
\item[(iii)] $a^\circ \cap bR=\{0\}$, $~^\circ a\cap Rc=\{0\}$, $R=abR+ c^\circ$ and $R=Rca+~^\circ b$;
 \item[(iv)]\emph{\cite[Proposition 2.7]{D}}$R=abR\oplus c^\circ$ and $R=Rca\oplus ~^\circ b$;
 \item[(v)]\emph{\cite[Proposition 2.7]{D}} $R=abR+ c^\circ$ and $R=Rca+~^\circ b$.
\end{enumerate}
\end{cor}

Before we present our next result, we first state the following auxiliary lemma.

\begin{lem}\emph{\cite{Han-Chen-92, Zhu-Zhang-Chen-LAA15}}\label{ke-onesideinverse-lem1}
Let $R$ be a $*$-ring with 1 and let $a\in R$. Then
\begin{enumerate}
\item[(i)] $a\{1, 3\} \neq \emptyset \Leftrightarrow a^*R=a^*aR \Leftrightarrow Ra=Ra^*a \Leftrightarrow R=Ra^*+~^\circ a \Leftrightarrow R=aR+(a^*)^\circ.$
\item[(ii)] $a\{1, 4\}\neq \emptyset \Leftrightarrow aR=aa^*R \Leftrightarrow Ra^*=Raa^* \Leftrightarrow R=Ra+~^\circ (a^*) \Leftrightarrow R=a^*R+a^\circ.$
\end{enumerate}
\end{lem}

Using Theorem \ref{ke-onesideinverse-thm1} and  Lemma \ref{ke-onesideinverse-lem1}, we can get the following result easily, which derives the relations between left, right $(b, c)$-inverse and $\{1, 3\}$, $\{1, 4\}$ and Moore-Penrose inverse.

\begin{prop}\label{ke-onesideinverse-prop3-mp}
Let $R$ be a $*$-ring with 1 and let $a\in R$. Then:
\begin{enumerate}
\item[(i)] $a$ is left $(a^*, 1)$-invertible if and only if $a$ has a $\{1, 4\}$-inverse; 
\item[(ii)] $a$ is right $(1, a^*)$-invertible if and only if $a$ has a $\{1, 3\}$-inverse; 
\item[(iii)] $a$ is Moore-Penrose invertible if and only if $a$ is left $(a^*, 1)$-invertible and right $(1, a^*)$-invertible.
\end{enumerate}
\end{prop}
\proof  (i). Using Theorem \ref{ke-onesideinverse-thm1}, $a$ is left $(a^*, 1)$-invertible if and only if $Ra^*=Raa^*$. By Lemma \ref{ke-onesideinverse-lem1}, this is equivalent to $a$ is $\{1, 4\}$-invertible.

       (ii). Similar discuss as (i).

       (iii). It is well known that $a\in R$ is Moore-Penrose invertible if and only if $a\in aa^*R\cap Ra^*a$ if and only if it is both $\{1, 3\}$ and $\{1, 4\}$-invertible. Using (i) and (ii), (iii) holds.
\qed



Now we consider the relations between left, right $(b, c)$-invertible and left, right $\pi$-regular and strongly $\pi$-regular elements.

\begin{prop}\label{ke-onesideinverse-prop3-drazin}
Let $R$ be a ring with 1 and let $a\in R$. Then for some $n\in \mathbb{N}$,
\begin{enumerate}
\item[(i)] $a$ is left $(a^n, 1)$-invertible if and only if $a$ is left $\pi$-regular; 
\item[(ii)] $a$ is right $(1, a^n)$-invertible if and only if $a$ is right $\pi$-regular; 
\item[(iii)] $a$ is strongly $\pi$-regular (i.e. Drazin invertible) if and only if $a$ is left $(a^n, 1)$-invertible and right $(1, a^n)$-invertible.
\end{enumerate}
\end{prop}
\proof (i). Using Theorem \ref{ke-onesideinverse-thm1},
       $a$ is left $(a^n, 1)$-invertible if and only if $Ra=Ra^n$.
By the definition of left $\pi$-regular in \cite{Azumaya-1954}, (i) holds.

     (ii). Analogously as (i).

     (iii). Combining (i) and (ii), we get (iii).
\qed

\vspace{0.2cm}
Note that Proposition \ref{ke-onesideinverse-prop3-drazin} provides the relations between left, right $(b, c)$-invertible and left, right regular and strongly regular (i.e. group invertible) elements, for $n=1$.

As an application of Theorems \ref{ke-onesideinverse-thm1}, some other special cases are as follows.

\begin{prop}\label{ke-onesideinverse-prop4}
Let $R$ be a $*$-ring with 1 and let $a$ be an element of $R$. Then
\begin{enumerate}
\item[(i)] $a$ is right $(a^*, 1)$-invertible if and only if $R=aa^*R$
 ($a$ is left $(1, a^*)$-invertible if and only if $R=Ra^*a$, resp.).
\item[(ii)] $a$ is left $(a, a^*)$-invertible if and only if $Ra=Ra^*a^2$ ($a$ is right $(a, a^*)$-invertible if and only if $a^*R=a^*a^2R$, resp.).
\item[(iii)] $a$ is left $(a^*, a)$-invertible if and only if $Ra^*=Ra^2a^*$ ($a$ is right $(a^*, a)$-invertible if and only if $aR=a^2a^*R$. resp.).
\end{enumerate}
\end{prop}

\begin{rem} \rm{Let  $a, b, c\in R$.
\begin{enumerate}
\item[(i)] According to \cite[Proposition 3.3]{W-C-C}, if $a$ is $(b, c)$-invertible, then $b, c\in R^-$. However, if $a$ is left (right, resp.) $(b, c)$-invertible or left (right, resp.) $(b, c)$-invertible, the condition $b, c\in R^-$ doesn't have to hold in general.
For example, let $R=\mathbb{Z}_8=\mathbb{Z}/8\mathbb{Z}$, $a=\overline{5}, b=\overline{0}, c=\overline{2}$. Then $\overline{5}$ is left $(\overline{0}, \overline{2})$-invertible and it is also right annihilator $(\overline{0}, \overline{2})$-invertible. But $\overline{2}$ is not regular element in $\mathbb{Z}_8$.
\item[(ii)]  In Proposition \ref{ke-onesideinverse-prop3-mp}, \ref{ke-onesideinverse-prop3-drazin} and \ref{ke-onesideinverse-prop4} we have the several choices for $b$ and $c$: 1, $a^*$, $a^n$, for $n\in \mathbb{N}$. However, in each of these propositions the case when $b\neq c$ is studied. For the above mentioned choices of $b$ and $c$, in the case when $b=c$, we have the following statements: \\
     (a) $a$ is left (right, resp.) $(1, 1)$-invertible if and only if it is left (right, resp.) invertible.\\
     (b) $a$ is left (right, resp.) $(a, a)$-invertible if and only if it is left (right, resp.) regular.\\
     (c) $a$ is left (right, resp.) $(a^*, a^*)$-invertible if and only if it is left (right, resp.) $*$-regular if and only if $a$ is Moore-Penrose invertible.\\
     (d) $a$ is left (right, resp.) $(a^n, a^n)$-invertible if and only if it is left (right, resp.) $\pi$-regular, for some $n\in \mathbb{Z}$.
\item[(iii)] If $Rc \subseteq Rb \subseteq Rcab$ ($bR \subseteq cR \subseteq cabR$, resp.) (or $Rc=Rb$ ($bR=cR$, resp.)), then $a$ is left (right, resp.) $(b, c)$-invertible if and only if it is left (right, resp.) invertible along $b$ ($c$, resp.).
\end{enumerate}}
\end{rem}

\section{One-sided $(b, c)$-inverse of a product of three elements}
In this section, we present several necessary and sufficient conditions for the existence of one-sided $(b, c)$-inverse of a product $paq$ in a ring $R$, where $p, a, q\in R$.

First, we consider the relation between the left $(b, c)$-inverse of $paq$ and one-sided generalized inverse of $pa$ and $aq$, where $p, a, q\in R$.

\begin{thm}\label{ke-onesideinverse-thm5}
Let  $a, b, c, p, q\in R$. Then the following are equivalent:
\begin{enumerate}
\item[(i)] $paq$ is left $(b, c)$-invertible;
\item[(ii)] $pa$ is left $(qb, c)$-invertible and $aq$  is left $(b, cp)$-invertible.
\end{enumerate}
Furthermore, if $y\in R$ is a left $(b, c)$-inverse of $paq$, $x\in R$ is a left  $(qb, c)$-inverse of $pa$ and $z\in R$ is a left $(b, cp)$-inverse of $aq$, then the following relations hold
$$y=zax, \qquad x=qy \qquad \mbox{and} \qquad z=yp.$$
\end{thm}
\proof
    (i) $\Rightarrow$ (ii). Let $y\in R$ be a left $(b, c)$-inverse of $paq$. By Definition \ref{ke-onesideinverse-def1}, we have $Ry\subseteq Rc$ and $ypaqb=b$. Let $x=qy$ and $z=yp$. Then we have
    $$Rx=Rqy \subseteq Ry\subseteq Rc \quad \mbox{and}\quad xpaqb=qb,$$
    $$Rz=Ryp\subseteq Rcp \quad \mbox{and}\quad  zaqb=b.$$

    (ii) $\Rightarrow$ (i). If (ii) holds, by Definition \ref{ke-onesideinverse-def1}, we have
    $$Rx\subseteq Rc,\quad xpaqb=qb, \quad Rz\subseteq Rcp,\quad zaqb=b.$$
    Let $y=zax$. Then $Ry=Rzax \subseteq Rx \subseteq Rc$ and $ypaqb=zaxpaqb=zaqb=b$.
\qed

\medskip

Similarly, we get the analogous result for the right $(b, c)$-inverse of $paq$.

\begin{thm}\label{ke-onesideinverse-thm5.1}
Let  $a, b, c, p, q\in R$. Then the following are equivalent:
\begin{enumerate}
\item[(i)] $paq$ is right $(b, c)$-invertible;
\item[(ii)] $pa$ is right $(qb, c)$-invertible and $aq$  is right $(b, cp)$-invertible.
\end{enumerate}
Moreover, if $y\in R$ is a right $(b, c)$-inverse of $paq$, $x\in R$ is a right  $(qb, c)$-inverse of $pa$ and $z\in R$ is a right $(b, cp)$-inverse of $aq$, then the following relations hold
$$y=zax, \qquad x=qy \qquad \mbox{and} \qquad z=yp.$$
\end{thm}

Note that Theorem \ref{ke-onesideinverse-thm5} (Theorem \ref{ke-onesideinverse-thm5.1}, resp.) provide necessary and sufficient conditions for the existence of left (right, resp.) inverse along element $d\in R$ of a product of three elements, in the case when $(b,c)=(d,d)$.



\medskip

In a similar way as in the above theorems, we obtain the following result.

\begin{thm}\label{ke-onesideinverse-thm6}
Let $a, b, c, p, q\in R$. If there exists $q'\in R$ such that $q'qb=b$, then the following are equivalent:
\begin{enumerate}
\item[(i)] $paq$ is left $(b, c)$-invertible;
\item[(ii)] $a$ is left $(qb, cp)$-invertible.
\end{enumerate}
Moreover, if $y\in R$ is a left $(b, c)$-inverse of $paq$ and $w\in R$ is a left $(qb, cp)$-inverse of $a$, then the following relation holds:
 $$w=qyp.$$
\end{thm}
\proof (i) $\Rightarrow$ (ii). If (i) holds, by Definition \ref{ke-onesideinverse-def1}, we have $Ry\subseteq Rc$ and $ypaqb=b$. Let $w=qyp$. Then we get $$Rw=Rqyp \subseteq Ryp\subseteq Rcp \quad\mbox{and}\quad waqb=(qyp)aqb=q(ypaqb)=qb.$$ This implies that $a$ has a left $(qb, cp)$-inverse $w$.

(ii) $\Rightarrow$ (i). Assume that $a$ has a left $(qb, cp)$-inverse $w$. By Definition \ref{ke-onesideinverse-def1}, we have $Rw\subseteq Rcp$ and $waqb=qb$. Then, there is $r\in R$ such that $w=rcp$. Since $q'qb=b$, multiplying by $q'$ on the left of $waqb=qb$ gives $b=q'waqb=q'(rcp)aqb$. Let $y=q'rc$. Then $Ry=Rq'rc\subseteq Rc$ and $b=(q'rc)paqb=ypaqb$. In Consequence, by Definition \ref{ke-onesideinverse-def1}, $paq$ has a left $(b, c)$-inverse $y$.
\qed


\medskip
Analogously, we can show the following characterization for the right $(b, c)$-inverse of $paq$.

\begin{thm}\label{ke-onesideinverse-thm6.1}
Let $a, b, c, p, q\in R$. If there exists $p'\in R$ such that $cpp'=c$, then the following are equivalent:
\begin{enumerate}
\item[(i)] $paq$ is right $(b, c)$-invertible;
\item[(ii)] $a$ is right $(qb, cp)$-invertible.
\end{enumerate}
Furthermore, if $y\in R$ is a right $(b, c)$-inverse of $paq$ and $w\in R$ is a right $(qb, cp)$-inverse of $a$, then the following relation holds:
$$w=qyp.$$
\end{thm}


Applying Theorem \ref{ke-onesideinverse-thm6}, Theorem \ref{ke-onesideinverse-thm6.1} and Corollary \ref{ke-onesideinverse-cor3}, we have the next corollary for the $(b, c)$-inverse of $paq$.

\begin{cor}\label{ke-onesideinverse-cor6.1}
Let $a, b, c, p, q\in R$. If there exist $p', q'\in R$ such that $q'qb=b$ and $cpp'=c$, then the following are equivalent:
\begin{enumerate}
\item[(i)] $paq$ is  $(b, c)$-invertible;
\item[(ii)] $a$ is $(qb, cp)$-invertible.
\end{enumerate}
In this case, if $y\in R$ is the  $(b, c)$-inverse of $paq$ and $w\in R$ is the $(qb, cp)$-inverse of $a$, then the following relation holds:$$w=qyp.$$
\end{cor}

Note that Theorem \ref{ke-onesideinverse-thm6} (Theorem \ref{ke-onesideinverse-thm6.1}, resp.) is also valid if we replace the condition ``there exists $q'\in R$ such that $q'qb=b$" (``there exists $p'\in R$ such that $cpp'=c$", resp.) with stronger condition ``$b$ is left invertible" (``$c$ is right invertible"). Hence, Corollary \ref{ke-onesideinverse-cor6.1} is valid if the condition ``$b$ is left invertible and $c$ is right invertible" holds instead of  ``there exist $p', q'\in R$ such that $q'qb=b$ and $cpp'=c$".

\medskip

In the following result we obtain the relation between the $(b, c)$-inverse of $paq$ and one-sided generalized inverses of $pa$ and $aq$.

\begin{thm}\label{ke-onesideinverse-thm7}
Let $a, b, c, p, q\in R$. If there exist $p', q'\in R$ such that $q'qc=c$ and $bpp'=b$, then the following are equivalent:
\begin{enumerate}
\item[(i)] $paq$ is $(b, c)$-invertible;
\item[(ii)] $pa$ is right $(qb, qc)$-invertible and $aq$ is left $(bp, cp)$-invertible.
\end{enumerate}
Moreover, if $y$ is the $(b, c)$-inverse of $paq$, $x$ is a right $(qb, qc)$-inverse of $pa$ and $z$ is a left $(bp, cp)$-inverse of $aq$, then the following relation holds:  $$y=zax.$$
\end{thm}
\proof  (i) $\Rightarrow$ (ii). Suppose that $paq$ is  $(b, c)$-invertible, by \cite[Theorem 2.2]{D}, we have $Rb=Rcpaqb$ and $cR=cpaqbR$, which imply $Rbp=Rcpaqbp$ and $qcR=qcpaqbR$. According to Theorem \ref{ke-onesideinverse-thm1}, we obtain that $pa$ is right $(qb, qc)$-invertible and $aq$ is left $(bp, cp)$-invertible.

(ii) $\Rightarrow$ (i). 
 If (ii) holds, using Definition \ref{ke-onesideinverse-def1}, we have
\begin{equation}\label{ke-onesideinverse-eq6}
xR\subseteq qbR, \quad qcpax=qc, \quad Rz\subseteq Rcp,\quad zaqbp=bp.
\end{equation}
Since  $q'qc=c$ and $bpp'=b$,  from (\ref{ke-onesideinverse-eq6}) we get:
\begin{equation}\label{ke-onesideinverse-eq7}
cpax=c,\quad zaqb=b.
\end{equation}
The condition $xR\subseteq qbR$ and $Rz\subseteq Rcp$ implies $x=qbr_1$ and $z=r_2cp$, for some $r_1, r_2\in R$.
Let $y=zax$. Then
$$
Ry=Rzax\overset{\tiny{(\ref{ke-onesideinverse-eq6})}}\subseteq Rcpax\overset{\tiny{(\ref{ke-onesideinverse-eq7})}}=Rc, \quad yR=zaxR\overset{\tiny{(\ref{ke-onesideinverse-eq6})}}\subseteq zaqbR\overset{\tiny{(\ref{ke-onesideinverse-eq7})}}=bR,
$$
$$ypaqb=(zax)paqb=(r_2cp)axpaqb\overset{\tiny{(\ref{ke-onesideinverse-eq7})}}
=r_2cpaqb=zaqb\overset{\tiny{(\ref{ke-onesideinverse-eq7})}}=b,$$
$$cpaqy=cpaq(zax)=cpaqza(qbr_1)\overset{\tiny{(\ref{ke-onesideinverse-eq7})}}=
cpaqbr_1=cpax\overset{\tiny{(\ref{ke-onesideinverse-eq7})}}=c.$$
Therefore, by Definition \ref{ke-onesideinverse-def1}, $paq$ has a left and right $(b, c)$-inverse $y$. Hence, by Corollary \ref{ke-onesideinverse-cor3}, $paq$ has a  $(b, c)$-inverse $y$.
\qed

\medskip

The above theorem is also valid if the condition "$b$ is right and $c$ is left invertible" holds instead of "there exist $p', q'\in R$ such that $q'qc=c$ and $bpp'=b$".
Note that we can get the related results for one-sided invertibility along an element, as a direct application of the above theorems.



\section{Some properties of one-sided $(b, c)$-inverses}

In this section, we investigate some properties of one-sided $(b, c)$-inverse in rings.

First, we need  the following auxiliary lemma.

\begin{lem}\emph{\cite[Exercise 1.6]{Lam-book01}}\label{ke-lem4.2}
Let $a, b \in R$. Then
\begin{enumerate}
\item[(i)] $1+ab$ is left invertible if and only if $1+ba$ is left invertible. Moreover, if $y(1 + ab) = 1$, then $(1 - bya)(1 + ba) = 1$.
\item[(ii)]  $1+ab$ is right invertible if and only if $1+ba$ is right invertible. Moreover, if $(1 + ab)x = 1$, then $(1 + ba)(1 - bxa) = 1$.
\item[(iii)] $1+ab$ is  invertible if and only if $1+ba$ is  invertible. Moreover, $(1 +ba)^{-1}=1-b(1 +ab)^{-1}a$.
\end{enumerate}
\end{lem}

Note that Lemma \ref{ke-lem4.2} (iii) is known as the Jacobson's Lemma.

Now we will consider the relation between left $(b, c)$-inverse of $\alpha\in R$ and classical one-sided inverses of $1+(\alpha-a)a^\otimes$ and $1+a^\otimes (\alpha-a)$.

\begin{thm}\label{ke-onesideinverse-thm8}
Let $a, b, c, \alpha\in R$ be such that $a$ has a $(b, c)$-inverse $a^\otimes$. 
 The following are equivalent:
\begin{enumerate}
     \item[(i)] $\alpha$ is left $(b, c)$-invertible;
     \item[(ii)] $\alpha$ is right annihilator $(b, c)$-invertible;
     \item[(iii)] $1+(\alpha-a)a^\otimes$ is left invertible;
     \item[(iv)] $1+a^\otimes (\alpha-a)$ is left invertible.
   \end{enumerate}
\end{thm}
\proof  (i) $\Leftrightarrow$ (ii). Since $a$ has a $(b, c)$-inverse $a^\otimes$, by \cite[Proposition 3.3]{W-C-C}, then $b, c\in R^-$. By Proposition \ref{ke-onesideinverse-prop0.1}, (i) is equivalent to (ii).

 (iii) $\Leftrightarrow$ (iv). Let $u=1+(\alpha-a)a^\otimes$ and $v=1+a^\otimes (\alpha-a)$. Applying Lemma \ref{ke-lem4.2}, the left invertibility of $u$ is equivalent to the left invertibility of $v$.

(i) $\Rightarrow$ (iii). By Theorem \ref{ke-onesideinverse-thm3}, $\alpha$ is left $(b, c)$-invertible if and only if $R=Rc\alpha +~^\circ b$. 
Since $a$ has a $(b, c)$-inverse $a^\otimes$, by \cite[Proposition 6.1]{D}, we know that $a^\otimes a a^\otimes=a^\otimes$, $bR=a^\otimes R$ and $Rc=Ra^\otimes$. Then $R=Rc\alpha + ~^\circ b=Ra^\otimes \alpha + ~^\circ(a^\otimes)$. For any $z\in R$, we have $z=z_1a^\otimes \alpha+z_2$, where $z_1\in Ra^\otimes a$ and $z_2\in ~^\circ(a^\otimes)$. Note that $z_1\in Ra^\otimes a$ implies $z_1=z_1a^\otimes a$. Let $t=z_1+z_2$. Then
    $$t(1+a^\otimes (\alpha-a))=(z_1+z_2)(1+a^\otimes (\alpha-a))=z_1+z_1a^\otimes (\alpha-a)+z_2=z_1+z_1a^\otimes(\alpha-a)+z_2=z.$$
    As $z\in R$ is arbitrary, let $z=1$. Then $1+a^\otimes (\alpha-a)$ is left invertible. By Lemma \ref{ke-lem4.2}, $1+(\alpha-a)a^\otimes$ is left invertible.

(iv) $\Rightarrow$ (i). If $v=1+a^\otimes (\alpha-a)$ is left invertible, then there exists $t\in R$ such that  $t(1+a^\otimes (\alpha-a))=1$. Since $a$ has a $(b, c)$-inverse $a^\otimes$, we have $a^\otimes a a^\otimes=a^\otimes$, $bR=a^\otimes R$ and $Rc=Ra^\otimes$, moreover, $b\in R^-$. Hence the condition $bR=a^\otimes R$ gives $~^\circ b=~^\circ (a^\otimes)$. Thus, $1-a^\otimes a\in ~^\circ (a^\otimes)=~^\circ b$. Therefore, for any $s\in R$,  $$s=st(1+a^\otimes (\alpha-a))=sta^\otimes \alpha +st(1-a^\otimes a)\in Ra^\otimes \alpha + ~^\circ b=Rc\alpha + ~^\circ b.$$
Consequently, $R= Rc\alpha + ~^\circ b$.  By Theorem \ref{ke-onesideinverse-thm3}, $\alpha$ is left $(b, c)$-invertible.
\qed

\medskip

Dually, we have the similar property for right $(b, c)$-inverse of $\alpha\in R$.
\begin{thm}\label{ke-onesideinverse-thm8.1}
Let $a, b, c, \alpha\in R$ be such that $a$ has a $(b, c)$-inverse $a^\otimes$. 
The following are equivalent:
\begin{enumerate}
     \item[(i)] $\alpha$ is right $(b, c)$-invertible;
     \item[(ii)] $\alpha$ is left annihilator $(b, c)$-invertible;
     \item[(iii)] $1+(\alpha-a)a^\otimes$ is right invertible;
     \item[(iv)] $1+a^\otimes (\alpha-a)$ is right invertible.
   \end{enumerate}
\end{thm}



Note that we can get the related results for one-sided invertibility along an element.

\vspace{0.2cm} \noindent {\large\bf Acknowledgments}

           The first author is grateful to China Scholarship Council for supporting her to purse her further study with Professor D. S. Cvetkovi\'{c}-Ili\'{c} in University of Ni\v{s}, Serbia. 
           The research was supported by The National Natural Science Foundation of China (No. 11371089),
           the Specialized Research Fund for the Doctoral Program of Higher Education (No. 20120092110020),
           the Natural Science Foundation of Jiangsu Province (No. BK20141327),
           the Foundation of Graduate Innovation Program of Jiangsu Province (No. KYLX$_{-}$0080), the Ministry of Education, Sciences and Technological Development, Republic of Serbia (Grant No. 174007).


\begin{thebibliography}{99}
\bibitem{Azumaya-1954}G. Azumaya, Strongly $\pi$-regular rings, J. Fac. Sci. Hokkaido Univ. 13 (1954) 34-39.




\bibitem{Benitez-Boasso-1507}J. Benitez, E. Boasso, The inverse along an element in rings, arXiv: 1507.05410 (2015).

\bibitem{Benitez-Boasso-1509}J. Benitez,  E. Boasso, The inverse along an element in rings with an involution, Banach algebras and $C^*$-algebras,  arXiv: 1509.04251 (2015).




\bibitem{Campbell-Meyer-2009} S. L. Campbell, C. D. Meyer, Generalized inverses of linear transformations, Philadelphia, SIAM, 2009.



\bibitem{Han-Chen-92}R. Z. Han, J. L. Chen, Generalized inverses of matrices over rings, Chinese Quarterly J. Math. 7 (4) (1992) 40-49.



\bibitem{Drazin1958} M. P. Drazin, Pseudo-inverses in associative rings and semigroups, Amer. Math. Monthly 65 (1958) 506-514.

\bibitem{D} M. P. Drazin, A class of outer generalized inverses, Linear Algebra Appl. 436 (2012) 1909-1923.

\bibitem{Drazin-LMA2014} M. P. Drazin, Commuting properties of generalized inverses, Linear Multilinear Algebra, 61 (12) (2013) 1675-1681.

\bibitem{Drazin-LAA14-Uniqueproof} M. P. Drazin, Generalized inverses: Uniqueness proofs and three new classes, Linear Algebra Appl. 449 (2014) 402-416.




\bibitem{Kantun-Montiel-LMA14}G. Kant\'{u}n-Montiel, Outer generalized inverses with prescribed ideals, Linear Multilinear Algebra 62 (9)(2014) 1187-1196.



\bibitem{Hartwig-1976block} R. E. Hartwig, Block generalized inverses, Arch. Ration. Mech. Anal. 61 (1976) 197-251.

\bibitem{Ke-Filomat} Y. Y. Ke, Y. F. Gao, J. L. Chen,  Representations of $(b, c)$-inverses in rings with involution, Filomat, (accepted).




\bibitem{Lam-book01}T. Y. Lam, A First Course in Noncommutative Rings, Second ed., Grad. Text in Math., Vol. 131, Springer-Verlag, Berlin-Heidelberg-New York, 2001.

\bibitem{Mary2011} X. Mary, On generalized inverses and Green's relations, Linear Algebra Appl. 434 (8) (2011) 1836-1844.

\bibitem{Mary2012}X. Mary, Natural generalized inverse and core of an element in semigroup, rings and Banach and operator algebras, Eur. J. Pure Appl. Math. 5 (2012), 160-173.

\bibitem{Mary-P 2012} X. Mary, P. Patr\'{i}cio, The inverse along a lower triangular matrix,  Appl. Math. Comput. 219 (2012) 886-891.

\bibitem{Mary-P 2013} X. Mary, P. Patr\'{i}cio, Generalized inverses modulo $\mathcal{H}$ in semigroups and rings, Linear Multilinear Algebra 61 (8) (2013) 1130-1135.



\bibitem{P} R. Penrose, A generalized inverse for matrices, Math. Proc. Cambridge Philos. Soc. 51 (1955) 406-413.








\bibitem{W-C-C} L. Wang, J. L. Chen, N. Castro-Gonz\'{a}lez, Characterizations of the $(b, c)$-inverse in a ring, submitted.



\bibitem{Zhu-Zhang-Chen-LAA15} H. H. Zhu, X. X. Zhang,  J. L. Chen, Generalized inverses of a facterization in a ring with involution, Linear Algebra Appl. 472 (2015) 142-150.

\bibitem{Zhu-Chen-Patricio-LMA15104}H. H. Zhu, J. L. Chen, P. Patr\'{i}cio, Further results on the inverse along an element in
semigroups and rings, Linear Multilinear Algebra, 64 (3) (2016) 393-403.

\bibitem{Zhu-Patricio-Chen-Zhang-LMA15105}H. H. Zhu, P. Patr\'{i}cio, J. L. Chen, Y. L. Zhang, The inverse along a product and its
applications, Linear Multilinear Algebra 64 (5) (2016) 834-841.

\end{thebibliography}
\end{document}